\documentclass[11pt]{article}
\usepackage{amsfonts}
\usepackage{mathrsfs}

\usepackage[latin1]{inputenc}
\usepackage{amsmath,amssymb}
\usepackage{latexsym}
\usepackage[active]{srcltx}
\usepackage[
bookmarks=true,         %generates bookmarks for all entries in the "Table of Contents"
bookmarksnumbered=true, %includes chapter-/header-/section-/subsection-/... -
colorlinks=true, pdfstartview=FitV, linkcolor=blue, citecolor=blue,
urlcolor=blue]{hyperref}

\newtheorem{theorem}{Theorem}[section]

\newtheorem{remark}[theorem]{Remark}
\numberwithin{equation}{section}
\def\sgn{{\mathop {{\rm sgn\, }}}}

\def\square{{\vcenter{\vbox{\hrule height.3pt
        \hbox{\vrule width.3pt height5pt \kern5pt
           \vrule width.3pt}
        \hrule height.3pt}}}}

\def\tlint{{- \kern-0.85em \int \kern-0.2em}}
\def\dlint{{- \kern-1.05em \int \kern-0.4em}}

\def \eref#1{\hbox{(\ref{#1})}}

\def \eref#1{\hbox{(\ref{#1})}}

\newenvironment{proof}[1][Proof]{\noindent\textit{#1.} }{\hfill \rule{0.5em}{0.5em}}

\begin{document}

\title{A limit theorem for some linear processes with innovations in the domain of attraction of a stable law}
\date{\today}
\author{Fangjun Xu\thanks{F. Xu is partially supported by National Natural Science Foundation of China (Grant No.11871219).} \\
}
\maketitle
\begin{abstract} Let $X=\{X_n: n\in\mathbb{N}\}$ be a linear process in which the coefficients are of the form $a_i=i^{-1}\ell(i)$ with $\ell$ being a slowly varying function at the infinity and the innovations are independent and identically distributed random variables belonging to the domain of attraction of an $\alpha$-stable law with $\alpha\in (1, 2]$. We will establish the asymptotic behavior of the partial sum process
\[
\bigg\{\sum\limits_{n=1}^{[Nt]} X_n: t\geq 0\bigg\}
\]
as $N$ tends to infinity, where $[t]$ is the integer part of the non-negative number $t$.

\noindent


\end{abstract}

\section{Introduction}

Let $X=\{X_n: n\in\mathbb{N}\}$ be a linear process defined by
\begin{align} \label{lp}
X_n=\sum^{\infty}_{i=1} a_i\varepsilon_{n-i}, 
\end{align}
where the innovations $\varepsilon_i$ are independent and identical distributed (i.i.d.) random variables belonging to the domain of attraction of an $\alpha$-stable law with $\alpha\in(1,2]$, $\mathbb{E}\, \varepsilon_1=0$ and $\mathbb{E}\, \varepsilon^2_1=\infty$, $a_i=i^{-1} \ell(i)$ with $\ell$ being a slowly varying function at infinity.  For the innovation $\varepsilon_1$, by Theorems 2.6.1 and 2.6.2 in \cite{IL}, there exist nonnegative constants $\sigma_1$ and $\sigma_2$ such that
\begin{align} \label{tail}
\mathbb{P}(\varepsilon_1\leq -x)&=(\sigma_1+o(1))x^{-\alpha} h(x)\quad \text{and}\quad \mathbb{P}(\varepsilon_1>x)=(\sigma_2 +o(1))x^{-\alpha} h(x)
\end{align}
as $x$ tends to infinity, where  $\sigma_1, \sigma_2\geq 0$, $\sigma_1+\sigma_2>0$ and $h(x)$ is a positive slowly varying function at infinity. For $\alpha\in (1,2]$,  it is well known that (\ref{tail}) is equivalent to 
\begin{align} \label{char}
\ln\mathbb{E} e^{\iota u \varepsilon_1}=-(\sigma+o(1))|u|^{\alpha}H(|u|^{-1})(1-\iota D \sgn(u))
\end{align}
as $u$ tends to $0$, where $\iota=\sqrt{-1}$, $\sigma=(\sigma_1+\sigma_2)\Gamma(|\alpha-1|)\cos(\frac{\pi\alpha}{2})$, $D=\frac{\sigma_1-\sigma_2}{\sigma_1+\sigma_2}\tan(\frac{\pi\alpha}{2})$ and 
\begin{align} \label{H}
H(t)=
\begin{cases}
h(t) & \text{for}\; \alpha\in (1,2),\\  \\
-\int^t_0 s^2 d\frac{h(s)}{s^2} & \text{for}\; \alpha=2
\end{cases}
\end{align}
(see, e.g., \cite{A}). Note that $H(\cdot)$ is a slowly varying function at infinity. So $\sum\limits^{\infty}_{i=1} |a_i|^{\alpha}H(|a_i|^{-1})<\infty$ for each $\alpha\in (1,2]$. Then, according to \cite{A} (or Proposition 5.4 in \cite{BJL}), the linear process $X=\{X_n: n\in\mathbb{N}\}$ defined in \eref{lp} is a.s. convergent.  By Definition 2.1 in \cite{SSX} and (\ref{char}),  $X=\{X_n: n\in\mathbb{N}\}$ is a long memory linear process when $\alpha\in (1,2)$. However, in the case $\alpha=2$, the memory property of $X=\{X_n: n\in\mathbb{N}\}$ is determined by the slowly varying functions $\ell(\cdot)$ and $H(\cdot)$. For example, if $H(x)\sim \ln x$ where $\sim$ indicates that the ratio of both sides tends to 1 as $x$ tends to infinity, then, by Definition 2.1 in \cite{SSX} and (\ref{char}), $X=\{X_n: n\in\mathbb{N}\}$ has long memory when $\ell(x)\sim 1$ and short memory when $\ell(x)\sim \ln^{-2} x$.  For simplicity of notation, we assume that $\ell(\cdot)$ is a positive slowly varying function at infinity.

There are a lot of papers on limit theorems for linear processes. However, only a small number of them are for linear processes with infinite variance innovations, in particular, when the innovations belong to the domain of attraction of an $\alpha$-stable law with $\alpha\in (0,2]$, see \cite{S, A, DR, MP,  BJL, PSXY}. The limit behavior of $\sum\limits_{i\in \mathbb{Z}} b_{n, i}\varepsilon_i$, under certain conditions on the coefficients $\{b_{n, i}: i\in\mathbb{Z}\}$, was studied in \cite{S}. Functional limit theorems for linear processes with absolutely summable and non-absolutely summable coefficients were given in \cite{A}. But specific forms of coefficients are required in \cite{A}.   A limit theorem for the partial sums of linear processes with $\alpha\in (0,2)$ and absolutely summable coefficients was established in \cite{DR}. Functional limit theorems for linear processes with absolutely summable coefficients were studied in \cite{BJL}. Limit theorems for linear random fields were studied in \cite{MP} and \cite{PSXY}.  Summable coefficients and $\alpha\in(1,2)$ are required in \cite{MP}. Both absolutely summable and non-absolutely summable coefficient cases are considered in the recent paper \cite{PSXY} under very mild conditions. However, the functional limit theorem for linear process $X=\{X_n: n\in\mathbb{N}\}$ defined in (\ref{lp}) is unknown.

In this paper, for the linear process $X=\{X_n: n\in\mathbb{N}\}$ defined in (\ref{lp}), we consider the asymptotic behavior of the partial sum process $\big\{\sum\limits_{n=1}^{[Nt]} X_n: t\geq 0\big\}$ as $N$ tends to infinity. Clearly, the linear process $X=\{X_n: n\in\mathbb{N}\}$ defined in (\ref{lp}) corresponds to the one with $\beta=1$ in \cite{A}. However, the functional limit theorem for the partial sum of our linear process is not considered in \cite{A}.  When $t=1$, the asymptotic behavior of the partial sum $\sum\limits_{n=1}^{N} X_n$ for our linear process $X=\{X_n: n\in\mathbb{N}\}$ would probably be covered by Theorem 3.3 in \cite{PSXY}. However, the concise expression of the normalizer is not easy to obtain, see the definition of the normalizer in (16) of \cite{PSXY}. We will give the concise expression of the normalizer and establish the convergence of finite-dimensional distributions.  Moreover, the expression of our normalizer in (\ref{normalizer}) implies that the techniques used in \cite{A} could not help us establish the desired functional limit theorem.  Therefore, for our linear process $X=\{X_n: n\in\mathbb{N}\}$, to establish the functional limit theorems for $\Big\{\sum\limits_{n=1}^{[Nt]} X_n: t\geq 0\Big\}$, we need new ideas.

To address the functional limit theorems for our linear process, we introduce some notation. Let $Z=\{Z_t: \; t\geq 0\}$ be an $\alpha$-stable Levy process with $\alpha\in (1,2]$ and 
\begin{align} \label{limit}
\mathbb{E}\, e^{\iota u Z_t}=\exp\Big(-t \sigma |u|^{\alpha}\big(1-\iota D \sgn(u)\big)\Big),
\end{align}
where $\sigma$ and $D$ are defined in (\ref{char}).

For the function $H$ defined in (\ref{H}), there exists a unique slowly varying function $H_{\alpha}$ such that
\begin{align} \label{Halpha}
H(N^{\frac{1}{\alpha}} H^{\frac{1}{\alpha}}_{\alpha}(N))\sim H_{\alpha}(N)
\end{align}
as $N$ tends to infinity (see, e.g., \cite{A}).

The following is the main result of this paper.

\begin{theorem} \label{thm1} Let $X=\{X_n: n\in\mathbb{N}\}$ be the linear process defined as in (\ref{lp}). Then
\[
\Big\{ A^{-1}_N \sum^{[Nt]}_{n=1} X_n:\; t\geq 0  \Big\} \overset{f.d.d.}{\longrightarrow} \{Z_t:\; t\geq 0\}
\]
as $N$ tends to infinity, where $\overset{f.d.d.}{\longrightarrow}$ denotes the convergence of finite dimensional distributions, 
\begin{align} \label{normalizer}
A_N=N^{\frac{1}{\alpha}}H^{\frac{1}{\alpha}}_{\alpha}(N) \sum^N_{i=1} i^{-1}\ell(i)
\end{align}
and $Z=\{Z_t:\; t\geq 0\}$ is the $\alpha$-stable Levy process defined in (\ref{limit}).
\end{theorem}

\begin{remark} 
The normalizer in (\ref{normalizer}) can be rewritten as 
\begin{align*} \label{normalizer}
A_N=N^{\frac{1}{\alpha}}H^{\frac{1}{\alpha}}_{\alpha}(N) \sum^N_{i=1} a_i.
\end{align*}
By Theorem 2.1(a) in \cite{R},  the above normalizer $A_N$ gives a universal form of the normalizers in Theorem 1 in \cite{A}. Moreover, a proper constant multiple of the above normalizer $A_N$ would probably be the concise expression of the normalizer in Theorem 3.3 of \cite{PSXY}.
\end{remark}

Throughout this paper, if not mentioned otherwise, the letter $c$ with or without a subscript, denotes a generic positive finite constant whose exact value is independent of $n$ and may change from line to line. We use $\iota$ to denote the imaginary unit $\sqrt{-1}$. For a complex number $z$, we use $|z|$ to denote its modulus. Moreover, we use $[x]$ denotes the integer part of $x\geq 0$.

\section{Proof of Theorem \ref{thm1}}
\begin{proof} For each $m\in\mathbb{N}$, fix some numbers $u_i\in\mathbb{R}$ and $t_i\geq 0$ for $i=1,\cdots,m$. For simplicity, we assume that $u_i\neq 0$ for each $i=1,\cdots,m$ and $0=t_0<t_1<t_2<\cdots<t_m<\infty$. Let 
\[
Y_N(t)=A^{-1}_N \sum\limits^{[Nt]}_{n=1} X_n
\] 
for each $t\geq 0$.

It is easy to see that
\begin{align*}
\mathbb{E} \exp\left(\iota \sum\limits^m_{i=1} u_i Y_N(t_i)\right)= \mathbb{E} \exp\left(\iota \sum\limits^m_{i=1} v_i \big(Y_N(t_i)-Y_N(t_{i-1})\big)\right)
\end{align*}
and 
\begin{align*}
\mathbb{E} \exp\left(\iota \sum\limits^m_{i=1} v_i \big(Z_{t_i}-Z_{t_{i-1}}\big)\right)=\exp\left(-\sum\limits^m_{i=1} (t_i-t_{i-1}) \sigma |v_i|^{\alpha}\big(1-\iota D \sgn(v_i)\big)\right),
\end{align*}
where $v_i=\sum\limits^m_{j=i} u_i$ for $i=1,2,\cdots, m$. 

So we only need to show that 
\begin{align} \label{convergence}
&\lim_{N\to\infty} \mathbb{E} \exp\left(\iota \sum\limits^m_{i=1} v_i \big(Y_N(t_i)-Y_N(t_{i-1})\big)\right) \nonumber\\
&\qquad\qquad\qquad\qquad=\exp\left(-\sum\limits^m_{i=1} (t_i-t_{i-1}) \sigma |v_i|^2\big(1-\iota D \sgn(v_i)\big)\right).
\end{align}

For each $i=1,2,\cdots, m$, define
\[
a^{[Nt_i]}_j=\sum\limits^{[Nt_i]}_{n=(j+1)\vee [Nt_{i-1}]+1}a_{n-j}.
\] 
Then
\begin{align*}
Y_N(t_i)-Y_N(t_{i-1})=A^{-1}_N\sum^{[Nt_i]-1}_{j=-\infty} a^{[Nt_i]}_j\varepsilon_j=A^{-1}_N\sum^{[Nt_i]-1}_{j=[Nt_{i-1}]} a^{[Nt_i]}_j\varepsilon_j+A^{-1}_N\sum^{[Nt_{i-1}]-1}_{j=-\infty} a^{[Nt_i]}_j\varepsilon_j
\end{align*}
and
\begin{align*}
\sum\limits^m_{i=1} v_i \Big(Y_N(t_i)-Y_N(t_{i-1})\Big)
&=\sum\limits^m_{i=1} \left[v_i A^{-1}_N \sum^{[Nt_i]-1}_{j=[Nt_{i-1}]} a^{[Nt_i]}_j\varepsilon_j+v_i A^{-1}_N  \sum^{[Nt_{i-1}]-1}_{j=-\infty} a^{[Nt_i]}_j\varepsilon_j \right]\\
&=\sum^{-1}_{j=-\infty} \left(\sum^m_{i=1} v_i A^{-1}_N a^{[Nt_i]}_j \right) \varepsilon_j +\sum^m_{k=1}\sum^{[Nt_k]-1}_{j=[Nt_{k-1}]} \left(\sum^m_{i=k} v_i A^{-1}_N a^{[Nt_i]}_j \right) \varepsilon_j.
\end{align*}

By Proposition 1(i) in \cite{A} and Theorem 2.6.1(a) in \cite{R}, we could easily obtain that 
\[
\lim\limits_{N\to\infty} A^{-1}_N \sup_{j} a^{[Nt_i]}_j=0
\] 
for each $i=1,2,\cdots,m$. Hence, for $N$ large enough, (\ref{char}) implies that
\begin{align*}
\ln \mathbb{E} \exp\left(\iota \sum\limits^m_{i=1} v_i A^{-1}_N\Big(Y_N(t_i)-Y_N(t_{i-1})\Big)\right)=(\sigma+o(1))(\text{I}_{N,1}+\text{I}_{N,2}),
\end{align*}
where
\begin{align*}
\text{I}_{N,1}=\sum\limits^{-1}_{j=-\infty} \Big|\sum^m_{i=1} v_i A^{-1}_N a^{[Nt_i]}_j \Big|^{\alpha}H\Big(\big|\sum^m_{i=1} v_i A^{-1}_N a^{[Nt_i]}_j \big|^{-1}\Big)\Big(1-\iota D \sgn\Big(\sum^m_{i=1} v_i A^{-1}_N a^{[Nt_i]}_j\Big)\Big)
\end{align*}
and
\begin{align*}
\text{I}_{N,2}=\sum^m_{k=1}\sum^{[Nt_k]-1}_{j=[Nt_{k-1}]}  \Big|\sum^m_{i=k} v_i A^{-1}_N a^{[Nt_i]}_j \Big|^{\alpha}H\Big(\big|\sum^m_{i=k} v_i A^{-1}_N a^{[Nt_i]}_j\big|^{-1}\Big)\Big(1-\iota D \sgn\Big(\sum^m_{i=k} v_i A^{-1}_N a^{[Nt_i]}_j\Big)\Big).
\end{align*}

In the sequel, we will show that
\begin{align} \label{in1}
\lim_{N\to\infty}\ \text{I}_{N,1}=0
\end{align}
and 
\begin{align} \label{in2}
\lim_{N\to\infty}\ \text{I}_{N,2}=\sum\limits^m_{k=1} (t_k-t_{k-1}) |v_k|^{\alpha}\big(1-\iota D \sgn(v_k)\big).
\end{align}

\noindent
{\bf Proof of (\ref{in1})}:  It suffices to show that
\begin{align} \label{k=0}
\lim_{N\to\infty}\sum\limits^{-1}_{j=-\infty} \Big|\sum^m_{i=1} v_i A^{-1}_N  a^{[Nt_i]}_j \Big|^{\alpha}H\Big(\big|\sum^m_{i=1} v_i A^{-1}_N a^{[Nt_i]}_j \big|^{-1}\Big)=0.
\end{align}
Note that
\begin{align} \label{iin12}
\sum\limits^{-1}_{j=-\infty} \Big|\sum^m_{i=1} v_i A^{-1}_N  a^{[Nt_i]}_j \Big|^{\alpha}H\Big(\big|\sum^m_{i=1} v_i A^{-1}_N a^{[Nt_i]}_j \big|^{-1}\Big)=\text{II}_{N,1}+\text{II}_{N,2},
\end{align}
where 
\begin{align*}
\text{II}_{N,1}=\sum\limits^{-N-1}_{j=-\infty} \Big|\sum^m_{i=1} v_i A^{-1}_N  a^{[Nt_i]}_j \Big|^{\alpha}H\Big(\big|\sum^m_{i=1} v_i A^{-1}_N a^{[Nt_i]}_j \big|^{-1}\Big)
\end{align*}
and 
\begin{align*}
\text{II}_{N,2}=\sum\limits^{-1}_{j=-N} \Big|\sum^m_{i=1} v_i A^{-1}_N  a^{[Nt_i]}_j \Big|^{\alpha}H\Big(\big|\sum^m_{i=1} v_i A^{-1}_N a^{[Nt_i]}_j \big|^{-1}\Big).
\end{align*}

Recall that $\lim\limits_{N\to\infty} \frac{H(N^{\frac{1}{\alpha}} H^{\frac{1}{\alpha}}_{\alpha}(N))}{H_{\alpha}(N)}=1$ and $\lim\limits_{N\to\infty} A^{-1}_N \sup_{j} a^{[Nt_i]}_j=0$. Then, for any $\eta\in(0,\alpha-1)$, 
\begin{align} \label{iin1}
\limsup_{N\to\infty}\text{II}_{N,1}
&=\limsup_{N\to\infty}\sum\limits^{-N-1}_{j=-\infty} N^{-1} \Bigg|\frac{\sum\limits^m_{i=1} v_i  a^{[Nt_i]}_j}{\sum\limits^N_{i=1} i^{-1}\ell(i)} \Bigg|^{\alpha}\frac{H\Bigg(N^{\frac{1}{\alpha}}H^{\frac{1}{\alpha}}_{\alpha}(N)\bigg| \frac{\sum\limits^N_{i=1} i^{-1}\ell(i)}{\sum\limits^m_{i=1} v_i a^{[Nt_i]}_j}\bigg|\Bigg)}{H\big(N^{\frac{1}{\alpha}}H^{\frac{1}{\alpha}}_{\alpha}(N)\big)}\frac{H\big(N^{\frac{1}{\alpha}}H^{\frac{1}{\alpha}}_{\alpha}(N)\big)}{H_{\alpha}(N)} \nonumber\\
&\leq c_1 \limsup_{N\to\infty} N^{-1}  \sum\limits^{-N-1}_{j=-\infty} \Bigg|\frac{\sum\limits^m_{i=1} a^{[Nt_i]}_j}{\sum\limits^N_{i=1} i^{-1}\ell(i)} \Bigg|^{\alpha-\eta}  \nonumber \\
%&\leq c_2  \limsup_{N\to\infty} N^{-1}  \sum\limits^{\infty}_{j=N+1}   \Bigg| \frac{\sum\limits^m_{i=1} a^{[Nt_i]}_j}{\sum\limits^N_{i=1} i^{-1}\ell(i)} \Bigg|^{\alpha-\eta}  \nonumber \\
&\leq c_2  \limsup_{N\to\infty} N^{-1}\left|\sum\limits^N_{i=1} i^{-1}\ell(i)\right|^{\eta-\alpha} \sum\limits^{\infty}_{j=N+1} \left| Nj^{-1}\ell(j) \right|^{\alpha-\eta}  \nonumber \\
&\leq c_3  \limsup_{N\to\infty} N^{-1}\left|\sum\limits^N_{i=1} i^{-1}\ell(i)\right|^{\eta-\alpha}   N\left|\ell(N) \right|^{\alpha-\eta}  \nonumber \\
&=c_3 \limsup_{N\to\infty}  \bigg|\frac{\ell(N)}{\sum\limits^N_{i=1} i^{-1}\ell(i)}\bigg|^{\alpha-\eta}=0,
\end{align}
where we use Proposition 2.6(ii) in \cite{R} in the first inequality, Theorem 2.1(a) in \cite{R} in the last inequality and Theorem 2.1(a) in \cite{R} in the last equality.

Moreover, for any $\delta\in (0,1)$ and $\eta\in(0,\alpha)$, we have
\begin{align}  \label{iin2}
\limsup_{N\to\infty}\text{II}_{N,2}
&=\limsup_{N\to\infty}\sum\limits^{-1}_{j=-N} N^{-1}H^{-1}_{\alpha}(N) \bigg| \frac{\sum\limits^m_{i=1} v_i  a^{[Nt_i]}_j}{\sum\limits^N_{i=1} i^{-1}\ell(i)} \bigg|^{\alpha}H\Bigg(N^{\frac{1}{\alpha}}H^{\frac{1}{\alpha}}_{\alpha}(N)\bigg|\frac{\sum\limits^N_{i=1} i^{-1}\ell(i)}{\sum\limits^m_{i=1}  v_i a^{[Nt_i]}_j}\bigg|\Bigg)  \nonumber \\
&\leq c_4 \limsup_{N\to\infty} N^{-1} \Bigg(\sum\limits^{-1}_{j=-[\delta N]} \bigg|\frac{\sum\limits^m_{i=1}  v_i  a^{[Nt_i]}_j}{\sum\limits^N_{i=1} i^{-1}\ell(i)} \bigg|^{\alpha-\eta}+ \sum\limits^{-[\delta N]-1}_{j=-N} \bigg|\frac{\sum\limits^m_{i=1}  v_i  a^{[Nt_i]}_j}{\sum\limits^N_{i=1} i^{-1}\ell(i)} \bigg|^{\alpha-\eta}\Bigg)  \nonumber \\
&\leq c_5 \limsup_{N\to\infty} N^{-1} \Bigg(\sum\limits^{-1}_{j=-[\delta N]} \bigg|\frac{\sum\limits^m_{i=1}  a^{[Nt_i]}_j}{\sum\limits^N_{i=1} i^{-1}\ell(i)} \bigg|^{\alpha-\eta}+ \sum\limits^{-[\delta N]-1}_{j=-N} \bigg|\frac{\sum\limits^m_{i=1}  a^{[Nt_i]}_j}{\sum\limits^N_{i=1} i^{-1}\ell(i)} \bigg|^{\alpha-\eta}\Bigg)  \nonumber \\
&\leq c_6\, \delta,
\end{align}
where in the first inequality we use Proposition 2.6(ii) in \cite{R} and in the last inequality we use the facts that 
\[
\sup\limits_{i=1,\cdots,m}\Bigg|\frac{a^{[Nt_i]}_j}{\sum\limits^N_{i=1} i^{-1}\ell(i)} \Bigg|\leq 1
\]
for $j=-[\delta N], \cdots, -1$ and 
\[
\lim\limits_{N\to\infty}\sup\limits_{i=1,\cdots,m}\frac{a^{[Nt_i]}_j}{\sum\limits^N_{i=1} i^{-1}\ell(i)}=0
\]
for $j=-N,\cdots, -[\delta N]-1$.
 
Combining (\ref{iin12}), (\ref{iin1}) and (\ref{iin2}) and letting $\delta\downarrow 0$, we could get the result in (\ref{k=0}) and thus the result in (\ref{in1}).

\medskip
\noindent
{\bf Proof of (\ref{in2})}: For each $k=1,\cdots, m$, we observe that
\begin{align} \label{iin34}
&\sum^{[Nt_k]-1}_{j=[Nt_{k-1}]}  \Big|\sum^m_{i=k} v_i  A^{-1}_N a^{[Nt_i]}_j \Big|^{\alpha}H\Big(\big|\sum^m_{i=k} v_i A^{-1}_N a^{[Nt_i]}_j\big|^{-1}\Big)\bigg(1-\iota D \sgn\Big(\sum^m_{i=k} v_i A^{-1}_N a^{[Nt_i]}_j\Big)\bigg) \nonumber\\
&\qquad\qquad\qquad:=\text{II}_{N,3}+\text{II}_{N,4},
\end{align}
where
\begin{align*}
\text{II}_{N,3}=\sum^{[\delta Nt_k]-1}_{j=[Nt_{k-1}]}  \Big|\sum^m_{i=k} v_i  A^{-1}_N a^{[Nt_i]}_j \Big|^{\alpha}H\Big(\big|\sum^m_{i=k} v_i A^{-1}_N a^{[Nt_i]}_j\big|^{-1}\Big)\Big(1-\iota D \sgn\Big(\sum^m_{i=k} v_i A^{-1}_N a^{[Nt_i]}_j\Big)\Big)
\end{align*}
and
\begin{align*}
\text{II}_{N,4}=\sum^{[Nt_k]-1}_{j=[\delta Nt_{k}]}  \Big|\sum^m_{i=k} v_i  A^{-1}_N a^{[Nt_i]}_j \Big|^{\alpha}H\Big(\big|\sum^m_{i=k} v_i A^{-1}_N a^{[Nt_i]}_j\big|^{-1}\Big)\Big(1-\iota D \sgn\Big(\sum^m_{i=k} v_i A^{-1}_N a^{[Nt_i]}_j\Big)\Big)
\end{align*}
with $\delta\in (0,1)$ and close enough to $1$.

We first estimate $\text{II}_{N,3}$. Here $j=[Nt_{k-1}],\cdots, [\delta N t_k]-1$. If $i=k$, then
\begin{align*}
\lim_{N\to\infty}\frac{a^{[Nt_k]}_j}{\sum\limits^N_{i=1} i^{-1}\ell(i)}
&=\lim_{N\to\infty}\frac{\sum\limits^{[Nt_k]}_{n=j+1} (n-j)^{-1}\ell(n-j)}{\sum\limits^N_{i=1} i^{-1}\ell(i)}=\lim_{N\to\infty}\frac{\sum\limits^{[Nt_k]-j}_{i=1} i^{-1}\ell(i)}{\sum\limits^N_{i=1} i^{-1}\ell(i)}=1,
\end{align*}
where in the last equality we use Theorem 2.1(a) in \cite{R} and the facts that 
\[
[Nt_k]-[\delta Nt_k]+1\leq [Nt_k]-j\leq [Nt_k]-[Nt_{k-1}],
\]
$\lim\limits_{N\to\infty}\frac{[Nt_k]-[\delta Nt_{k}]+1}{N}=t_k(1-\delta)>0$ and $\lim\limits_{N\to\infty}\frac{[Nt_k]-[Nt_{k-1}]}{N}=t_k-t_{k-1}>0$.

If $i=k+1,\cdots, m$, then
\begin{align*}
\lim_{N\to\infty}\frac{a^{[Nt_i]}_j}{\sum\limits^N_{i=1} i^{-1}\ell(i)}
&=\lim_{N\to\infty}\frac{\sum\limits^{[Nt_i]}_{n=[Nt_{i-1}]+1} (n-j)^{-1}\ell(n-j)}{\sum\limits^N_{i=1} i^{-1}\ell(i)}=\lim_{N\to\infty}\frac{\sum\limits^{[Nt_i]}_{i=1} i^{-1}\ell(i)-\sum\limits^{[Nt_{i-1}]}_{i=1} i^{-1}\ell(i)}{\sum\limits^N_{i=1} i^{-1}\ell(i)}\\
&=\lim_{N\to\infty}\frac{\sum\limits^{[Nt_i]}_{i=1} i^{-1}\ell(i)}{\sum\limits^N_{i=1} i^{-1}\ell(i)}-\lim_{N\to\infty}\frac{\sum\limits^{[Nt_{i-1}]}_{i=1} i^{-1}\ell(i)}{\sum\limits^N_{i=1} i^{-1}\ell(i)}\\
&=0,
\end{align*}
where in the last two equalities we use Theorem 2.1(a) in \cite{R} and the facts that $\lim\limits_{N\to\infty}\frac{[Nt_i]}{N}=t_i>0$ and $\lim\limits_{N\to\infty}\frac{[Nt_{i-1}]}{N}=t_{i-1}\geq t_k>0$. Therefore, 
\[
\lim\limits_{N\to\infty}\frac{\sum\limits^m_{i=k} v_i  a^{[Nt_i]}_j}{v_k  a^{[Nt_k]}_j}=1.
\]
 
Using properties of slowly varying functions, we can easily obtain 
\begin{align} \label{iin3}
&\lim\limits_{N\to\infty}\text{II}_{N,3}  \nonumber\\
&=\lim\limits_{N\to\infty} N^{-1} H^{-1}_{\alpha}(N)\sum^{[\delta Nt_k]-1}_{j=[Nt_{k-1}]}  \bigg| \frac{ \sum\limits^m_{i=k}v_i a^{[Nt_i]}_j }{\sum\limits^N_{i=1} i^{-1}\ell(i)}\bigg|^{\alpha}H\Bigg(N^{\frac{1}{\alpha}}H^{\frac{1}{\alpha}}_{\alpha}(N)\bigg|\frac{\sum\limits^m_{i=k}  v_i  a^{[Nt_i]}_j}{\sum\limits^N_{i=1} i^{-1}\ell(i)}\bigg|^{-1}\Bigg) \nonumber\\
&\qquad\qquad\qquad\qquad\times \bigg(1-\iota D \sgn\Big(\frac{\sum\limits^m_{i=k}  v_i a^{[Nt_i]}_j}{\sum\limits^N_{i=1} i^{-1}\ell(i)}\Big)\bigg) \nonumber \\
&=\lim\limits_{N\to\infty} N^{-1} H^{-1}_{\alpha}(N)\sum^{[\delta Nt_k]-1}_{j=[Nt_{k-1}]}  \bigg|\frac{v_k a^{[Nt_k]}_j }{\sum\limits^N_{i=1} i^{-1}\ell(i)}\bigg|^{\alpha}H\Bigg(N^{\frac{1}{\alpha}}H^{\frac{1}{\alpha}}_{\alpha}(N)\bigg|\frac{v_k  a^{[Nt_k]}_j}{\sum\limits^N_{i=1} i^{-1}\ell(i)}\bigg|^{-1}\Bigg) \bigg(1-\iota D \sgn\Big(\frac{v_k a^{[Nt_k]}_j}{\sum\limits^N_{i=1} i^{-1}\ell(i)}\Big)\bigg)\nonumber \\
&=(\delta  t_k-t_{k-1}) |v_k|^{\alpha}(1-\iota D \sgn(v_k)).
\end{align}

We next estimate $\text{II}_{N,4}$. Here $j=[\delta Nt_k],\cdots, [N t_k]-1$. If $i=k$, then
\begin{align*}
\limsup_{N\to\infty}\frac{a^{[Nt_k]}_j}{\sum\limits^N_{i=1} i^{-1}\ell(i)}
&=\limsup_{N\to\infty}\frac{\sum\limits^{[Nt_k]}_{n=j+1} (n-j)^{-1}\ell(n-j)}{\sum\limits^N_{i=1} i^{-1}\ell(i)}=\limsup_{N\to\infty}\frac{\sum\limits^{[Nt_k]-j}_{i=1} i^{-1}\ell(i)}{\sum\limits^N_{i=1} i^{-1}\ell(i)}\leq 1,
\end{align*}
where in the last inequality we use Theorem 2.1(a) in \cite{R}.

If $i=k+1,\cdots, m$, then
\begin{align*}
\lim_{N\to\infty}\frac{a^{[Nt_i]}_j}{\sum\limits^N_{i=1} i^{-1}\ell(i)}
&=\lim_{N\to\infty}\frac{\sum\limits^{[Nt_i]}_{n=[Nt_{i-1}]+1} (n-j)^{-1}\ell(n-j)}{\sum\limits^N_{i=1} i^{-1}\ell(i)}=\lim_{N\to\infty}\frac{\sum\limits^{[Nt_i]}_{i=1} i^{-1}\ell(i)-\sum\limits^{[Nt_{i-1}]}_{i=1} i^{-1}\ell(i)}{\sum\limits^N_{i=1} i^{-1}\ell(i)}=0.
\end{align*}

Therefore, for any $\eta\in (0,\alpha)$, by Proposition 2.6(ii) in \cite{R}, we have
\begin{align} \label{iin4}
&\limsup\limits_{N\to\infty} \left|\text{II}_{N,4}\right| \nonumber\\
&\leq 2\limsup\limits_{N\to\infty} N^{-1} \sum^{[Nt_k]-1}_{j=[\delta Nt_{k}]}  \bigg| \frac{\sum\limits^m_{i=k}v_i a^{[Nt_i]}_j }{\sum\limits^N_{i=1} i^{-1}\ell(i)}\bigg|^{\alpha}\frac{H\Bigg(N^{\frac{1}{\alpha}}H^{\frac{1}{\alpha}}_{\alpha}(N)\bigg| \frac{\sum\limits^m_{i=k} v_i  a^{[Nt_i]}_j}{\sum\limits^N_{i=1} i^{-1}\ell(i)}\bigg|^{-1}\Bigg)}{H\big(N^{\frac{1}{\alpha}}H^{\frac{1}{\alpha}}_{\alpha}(N)\big)}\frac{H\big(N^{\frac{1}{\alpha}}H^{\frac{1}{\alpha}}_{\alpha}(N)\big)}{H_{\alpha}(N)}\nonumber\\
&\leq c_7\, \limsup\limits_{N\to\infty} N^{-1}  \sum^{[Nt_k]-1}_{j=[\delta Nt_{k}]}  \bigg|\frac{v_k a^{[Nt_k]}_j }{\sum\limits^N_{i=1} i^{-1}\ell(i)}\bigg|^{\alpha-\eta} \nonumber\\
&\leq c_8\, (1-\delta) t_k.
\end{align}
Combining (\ref{iin34}), (\ref{iin3}), (\ref{iin4}) and letting $\delta\uparrow 1$, we could get the desired result in (\ref{in2}). This completes the proof. 
\end{proof}

\bigskip

$\begin{array}{cc}

\begin{minipage}[t]{1\textwidth}

{\bf Fangjun Xu}

\medskip

KLATASDS-MOE, School of Statistics, East China Normal University, Shanghai, 200062, China 

\medskip

NYU-ECNU Institute of Mathematical Sciences at NYU Shanghai, Shanghai, 200062, China

\medskip

\texttt{fangjunxu@gmail.com, fjxu@finance.ecnu.edu.cn}

\end{minipage}

\hfill

\end{array}$

\end{document}